\documentclass[11pt,a4paper, leqno]{article}
\usepackage{amsmath}
\usepackage{amsfonts}
\usepackage{mathtools}
\usepackage{amssymb}
\usepackage{comment}
\usepackage{pstricks,xy,pst-node}
\xyoption{all}
\usepackage{amsthm}
\usepackage{booktabs} 
\usepackage{mathrsfs} 
\usepackage{url}
\usepackage{caption} 
\usepackage{enumitem}
\usepackage[breaklinks,bookmarks=false, hidelinks,draft=false]{hyperref} 
\usepackage{libertine}
\usepackage{libertinust1math} 

\usepackage[T1]{fontenc}
\usepackage[english,italian,polish]{babel}
\usepackage[utf8]{inputenc}

\makeatletter
\let\c@table\c@figure 
\let\ftype@table\ftype@figure 
\makeatother

\numberwithin{equation}{section}

\newtheorem{theorem}[equation]{Theorem}

\newtheorem{proposition}[equation]{Proposition}

\newtheorem*{theorem*}{Theorem}

\theoremstyle{definition}
\newtheorem{definition}[equation]{Definition}
\newtheorem{remark}[equation]{Remark}
\newtheorem*{ack}{Acknowledgements}

\newtheorem*{question*}{Question}

\theoremstyle{remark}
\newtheorem{example}[equation]{Example}

\DeclareMathOperator{\hhh}{h}

\DeclareMathOperator{\Spec}{Spec}

\DeclareMathOperator{\Pic}{Pic}

\DeclareMathOperator{\HHH}{H}

\DeclareMathOperator{\coker}{coker}

\newcommand{\cH}{\mathcal{H}}
\newcommand{\cM}{\mathcal{M}}
\newcommand{\cP}{\mathcal{P}}
\newcommand{\cO}{\mathcal{O}}
\newcommand{\PP}{\mathbb{P}}
\newcommand{\cW}{\mathcal{W}}
\newcommand{\Mgn}{\cM_{g,n}}

\author{Hanieh Keneshlou\thanks{
\selectlanguage{polish}Mathematical Institute of the Polish Academy of Sciences, Jana i Jędrzeja Śniadeckich 8, 00-656 Warszawa, Poland}, Fabio Tanturri\thanks{
\selectlanguage{italian}Dipartimento di Matematica, Università di Genova, via Dodecaneso 35, 16146 Genova, Italy}}
\title{On the unirationality of moduli spaces of pointed curves}
\date{}

\begin{document}
\selectlanguage{english}
\maketitle
\begin{abstract}
We show that $\cM_{g,n}$, the moduli space of smooth curves of genus $g$ together with $n$ marked points, is unirational for $g=12$ and $2 \leq n\leq 4$ and for $g=13$ and $1 \leq n \leq 3$, by constructing suitable dominant families of projective curves in $\PP^1 \times \PP^2$ and $\PP^3$ respectively. We also exhibit several new unirationality results for moduli spaces of smooth curves of genus $g$ together with $n$ unordered points, establishing their unirationality for $g=11, n=7$ and $g=12, n =5,6$.
\end{abstract}
\section{Introduction}

The geometry of algebraic curves varying in families is a very fascinating and old topic, dating back to the nineteenth century. The interest around this subject naturally led to the definition of the moduli space $\cM_g$ of smooth curves of genus $g$ over the complex numbers. The study of its birational geometry (or the geometry of its Deligne--Mumford compactification $\overline{\cM}_g$) has become a very active research area, especially after the unexpected results of Harris--Mumford--Eisenbud \cite{HarrisMumfordKodaira,EisenbudHarrisKodaira}: they showed that $\overline{\cM}_g$ is of general type for $g \geq 24$, thus contradicting  a long-standing conjecture by Severi about its unirationality for any $g$. The unirationality for $g \leq 10$ being already implied by classical results by Severi \cite{SeveriClassification}, a great deal of work has ever since been devoted to the study of the birational geometry of $\overline{\cM}_g$ for the remaining cases, leading to the unirationality for $g \leq 14$ \cite{SernesiUnirationality, ChangRanUnirationality, VerraUnirationality} and the uniruledness for $g \leq 15$ \cite{ChangRanKodaira, BrunoVerraRationally}. In \cite{ChangRanKodaira} it was claimed that $\overline{\cM}_{16}$ is uniruled, but the proof has recently been shown flawed \cite{Tseng}; upper bounds have been newly proved for its Kodaira dimension in \cite{FarkasVerraKodaira16,AgostiniBarrosPencils}, showing in particular that it is not of general type. The cases $\overline{\cM}_{22}$ and $\overline{\cM}_{23}$ having being freshly shown to be of general type in \cite{FarkasJensenPayne}, nothing is known for the remaining genera $17 \leq g \leq 21$.

Along with $\cM_g$, the birational geometry of other moduli spaces has been considered and studied. For instance, one can consider isomorphism classes of curves together with additional structures such as finite maps to $\mathbb{P}^1$ (which leads to Hurwitz spaces, see, e.g., \cite{SchreyerTanturri, Mullane20}), or curves equipped with special line bundles. Being of independent interest, these spaces can also be used to shed further light on the geometry of the underlying spaces $\overline{\cM}_g$.

This paper concerns the study of the birational geometry of two of these spaces, namely the moduli space $\overline{\cM}_{g,n}$ parametrising stable $n$-pointed genus $g$ curves, and its quotient by the permutation group $S_n$, which is usually referred to as $\mathcal{C}_{g,n}$, the universal symmetric product of degree $n$, and parametrises stable curves of genus $g$ together with $n$ unordered points. The study of the birational properties of these moduli spaces can give extra tools for a better understanding of the properties of the spaces $\cM_g$. To highlight this worthiness and mention one application, we remark that the unirationality of $\cM_{14,2}$ was a key ingredient for proving that $\cM_{15}$ is rationally connected in \cite{BrunoVerraRationally}.

\subsection*{Main contributions}

The aim of this paper is to provide new unirationality results for the two aforementioned spaces in the range $11 \leq g \leq 13$.
Being interested only in their birational geometry, we will mostly deal with $\cM_{g,n}$ and $\cM_{g,n}^u$ instead of their compactifications. By convention, by the Kodaira dimensions of these spaces we will mean the Kodaira dimensions of their compactifications.

After the preliminary Section \ref{sec:preliminaries}, the first part of this paper (Section \ref{uniratM1213}) is devoted to proving new unirationality results for $\cM_{g,n}$. Since \cite{Logan}, it is known that for $g\geq 4$ the space $\cM_{g,n}$ is of general type for $n$ large enough. Thus, it is natural to try to characterise the birational geometry and to determine the Kodaira dimension of the finitely many remaining cases for each $g$. This problem has been investigated from many point of views and several contributions have been provided; we postpone to Section \ref{uniratM1213} a brief account on the known results.

Our main contribution is the following.

\begin{theorem}
\label{mainthmintro}
The space $\cM_{12,n}$ is unirational for $2 \le n\leq 4$; the space $\cM_{13,n}$ is unirational for $1 \leq n\leq 3$.
\end{theorem}

In particular, we prove of the unirationality of $\cM_{g,n}$ for genera $g$ and $n \geq 1$ for which the only known case was $\cM_{12,1}$ \cite{BallicoCasnatiFontanari}.

The approach we use is as follows. For a fixed genus $g\in \{12,13\}$, we produce a unirational family of projective curves of genus $g$ dominating $\cM_g$. Then we show that general elements of this family can be linked via hypersurfaces of a suitable degree to particular auxiliary curves. If wisely chosen, these auxiliary curves will be contained in more hypersurfaces of that degree, allowing us to reverse the process and impose a certain number of marked points on the elements of the original unirational family.

We use projective models in $\mathbb{P}^3$ and $\mathbb{P}^1 \times \mathbb{P}^2$ for $\cM_{13,n}$, $\cM_{12,n}$, respectively. The family of curves of genus $13$ is obtained building upon the unirationality of $\mathcal{C}_{10,6}$, granted by \cite{BarrosGeometry} and for which we exhibit an alternative proof in the second part of the paper. The construction of the family of curves of genus $12$ is based on a particular construction in \cite{GeissThesis}, which we recall in details. Two alternative proofs for the unirationality of $\cM_{12,n}$ for $n$ up to $3$, building upon a different, independent construction or on a family constructed in \cite{SernesiUnirationality}, are also provided.

Let us now consider $\mathcal{C}_{g,n}$. In this paper, we will denote it by $\cM_{g,n}^u$, to stress out that its elements correspond to the choice of $n$ \emph{unordered} points on a genus $g$ curve. In Section \ref{uniratmgnusection} we provide a brief account on the known results about its birational geometry. For what concerns the unirationality, it is known for $g\leq 9$; beside the results easily deducible from $\cM_{g,n}$, only the cases $g=10, n=6,7$ \cite{BarrosGeometry} are known.

Our main contribution in this framework is in Section \ref{uniratmgnusection}, where we prove the following

\begin{theorem}
\label{mgnuIntro}
The space $\cM_{g,n}^u$ is unirational for $g=11, n=7$, for $g=12$ and $2 \leq n \leq 6$, and for $g=13$ and $n \leq 3$.
\end{theorem}

Thanks to Theorem \ref{mainthmintro}, we need to discuss only the cases $g=11, n=7$ and $g=12, n=5,6$. For the case $\cM_{12,5}^u$ we adopt an ad hoc argument based on the unirationality of a component of the Hilbert scheme of curves of genus $9$ and degree $15$ in $\PP^6$. For the other two cases the unirationality of $\cM_{g,n}^u$ is achieved as follows. For a suitable choice of $m<n$, we exhibit a dominant unirational family of curves of genus $g$ and degree $2g-2-m$ in $\PP^{g-m-1}$, thus reproving the unirationality of $\cM_{g,m}^u$; we then show that by performing liaison forth and back we can impose a certain number $m'$ of additional points on these curves, yielding the unirationality of $\cM_{g,n}^u$ for $m < n \leq m+m'$.

With the same general approach, we also exhibit (Remark \ref{alternativeBarros}) a constructive alternative proof of the unirationality of $\cM_{g,n}^u$ for $g=10$ and $n=6,7$, already achieved in \cite{BarrosGeometry}.

\begin{ack}
The authors wish to thank the referee for many remarks and suggestions which have led to a global improvement of this paper. They also thank D.\ Agostini and I.\ Barros for pointing out relevant references.

The authors are grateful to the Max Planck Institute for Mathematics in the Sciences of Leipzig, Germany, where the initial phase of this work was carried out.
\end{ack}

\section{Preliminaries}
\label{sec:preliminaries}
 
For the sake of clearness, we quickly recall here a few facts about Brill--Noether theory, for which we refer to \cite{ACGH}, and a few facts on liaison of curves.
\subsection{Brill--Noether Theory}
\label{BNTheory}
Let $C$ denote a smooth curve of genus $ g $, and let $ d,r $ be non-negative integers. A linear series $g^r_d$ on $C$ of degree $d$ and dimension $r$ is a pair $(L,V)$, $L \in \Pic^d(C)$ being a line bundle of degree $d$ and $V \subseteq \HHH^0(C,L)$ an $(r+1)$-dimensional vector subspace of sections of $L$. A general curve $C$ of genus $g$ has a $g^r_d$ if and only if the Brill--Noether number
\[ \rho=\rho(g,r,d)=g-(r+1)(g+r-d) \]
is non-negative. In such case, the Brill--Noether scheme
\[ W^r_d(C)=\{ L \in \Pic^d(C) \mid \hhh^0(L) \ge r+1 \} \]
has dimension $\rho$. More generally, one can define the universal Brill--Noether scheme as
\[
\mathcal{W}^r_{g,d} = \{ (C,L) \mid C \in \mathcal{M}_{g}, L \in W^r_d(C) \}.
\]
\subsection{Liaison}
\begin{definition}
Let $ C $ and $ C^{\prime} $ be two curves in a projective variety $X$ of dimension $r$ with no embedded and no common components, contained in $r-1$ mutually independent hypersurfaces $Y_i \subset X$ which meet transversally. Let $Y$ denote the complete intersection curve $\cap Y_i$. The curves $ C $ and $ C^{\prime} $ are said to be \emph{geometrically linked} via $Y$ if $ C\cup C^{\prime}=Y$ scheme-theoretically.
\end{definition}
If $C$ and $C'$ are assumed to be locally complete intersections and to meet only in ordinary double points, then $\omega_Y|_C =\omega_C(C \cap C')$ and the arithmetic genera of the curves are related as follows:
\begin{equation*}
2(p_a(C)-p_a(C'))=\deg (\omega_C) - \deg(\omega_{C'})=\omega_X(Y_{1}+\cdots+Y_{r-1}).(C-C').
\end{equation*}
The above relation and the equality $\deg C + \deg C' = \deg Y$ can be used to deduce the genus and degree of $C'$ from the genus and degree of $C$.

\medskip
Let $X=\mathbb{P}^1 \times \mathbb{P}^2$ and $C$ be a curve of genus $p_a(C)$ and bidegree $(d_1, d_2)$. Let $Y_1, Y_2$ be two hypersurfaces of bidegree $(a_1,b_1)$ and $(a_2,b_2)$ containing $C$ and satisfying the above hypotheses.
The genus and the bidegree of $C'$ are 
\begin{equation}
\label{liaisonP1P2}
\begin{array}{rcl}
(d_1',d_2') &= &(b_1b_2-d_1,a_1b_2+a_2b_1-d_2),\\
p_a(C^{\prime}) &= &p_a(C) - \frac12 \left((a_1+a_2-2)(d_1-d_1')+(b_1+b_2-3)(d_2-d_2')\right).
\end{array}
\end{equation}

For curves embedded in a projective space $\mathbb{P}^r$, the invariants $p_a(C^{\prime}),d'$ of the curve $C'$ can be computed via
\begin{equation}
\label{liaisonPr}
\begin{array}{rcl}
d'&=&\prod d_i - d,\\
p_a(C^{\prime})&=&p_a(C)-\frac{1}{2}\left(\sum d_i-(r+1)\right)(d-d'),
\end{array}
\end{equation}
where the $d_i$'s are the degrees of the $r-1$ hypersurfaces $Y_i$ cutting out $Y$.

\subsection{Computational verifications}
In this paper we will often need to check on some explicit examples that some open conditions are generically satisfied. In order to do that, we will make use of the software \cite{M2}; the supporting code for this paper has been collected in \cite{pointedCurvesCode}, where an instance of its execution can be found. Although we could a priori run our computations directly on $\mathbb{Q}$, this can increase dramatically the required time of execution. Instead, we can work over a finite field $\mathbb{F}_p$, and view our choice of the initial parameters in $\mathbb{F}_p$ as the reduction modulo $p$ of some choices of parameters in $\mathbb{Z}$. The so-obtained example $E_p$ can be seen as the reduction modulo $p$ of a family of examples defined over an open part of $\Spec (\mathbb{Z})$. If $E_p$ satisfies an open condition, then a semicontinuity argument implies that the generic fiber $E$ satisfies the same open condition, and so does the general element of the family over $\mathbb{Q}$ or $\mathbb{C}$.

\section{On the unirationality of $\cM_{g,n}$}
\label{uniratM1213}

The moduli of curves $\cM_g$ are varieties of general type except for a finite number of cases, occurring for small values of $g$; the same principle holds true also for $\cM_{g,n}$, at least for $g>3$. Indeed, on the one hand Logan \cite{Logan} exhibited a natural number $\tau(g)$ for each fixed genus $g\geq 4$ such that $\cM_{g,n}$ is of general type for $n=\tau(g)$ (and hence, by the subadditivity of Kodaira dimensions, for $n \geq \tau(g)$). On the other hand, classical constructions of dominant families of curves allow us to prove the unirationality for small values of $g$ and $n$. 

In Table \ref{tableBiratGeom} we report on the known results: following an earlier notation, we have indicated for each genus $2 \leq g \leq 16$ the number $b(g)$ (respectively, $\sigma(g)$) such that $\cM_{g,n}$ is unirational (respectively, uniruled) for $n\leq b(g)$ (respectively, $n\leq \sigma(g)$). The number $\eta(g)$ (respectively, $\tau(g)$) denotes the smallest number of points $n$ for which the Kodaira dimension of $\cM_{g,n}$ is known to be non-negative (respectively, maximal). The table is based on contributions given by \cite{Logan, BrunoVerraRationally, VerraUnirationality,BiniFontanari, CasnatiFontanari, BallicoCasnatiFontanari,  BarrosGeometry, KeneshlouTanturri} for the unirationality; \cite{Logan, FarkasPopa, CasnatiFontanari, FarkasVerraUniversalJacobians,  Benzo, AgostiniBarrosPencils} for the uniruledness; \cite{Logan, FarkasKoszul,  FarkasVerraUniversalJacobians, FarkasVerraUniversalThetaDivisor,BarrosMullaneTwoModuli,SchwarzKodaira} for $\eta(g)$ and $\tau(g)$.

\begin{table}[h!bt]
	\begin{center}
		\setlength{\tabcolsep}{5pt}		
		\begin{tabular}{c|ccccccccccccccc} \toprule
			$g$ & 2 & 3& 4& 5& 6& 7& 8& 9& 10& 11& 12& 13& 14& 15& 16\\
			\midrule
			$b(g)$ & 12 & 14 & 15 & 12 & 15 & 11 & 11 & 9 & 5 & 6 & 1 & 0 & 2 \\
			$\sigma(g)$ & 13 & 14 & 15 & 14 & 15 & 13 & 12 & 10 & 9 & 10 & 7 & 4 & 3 & 2 \\
			$\eta(g)$ & 14 & & 16 & 15 & 16 & 14 & 14 & 13 & 10 & 11 & 10 & 11 & 10 & 10 & 9\\
			$\tau(g)$ & 15 & & 16 & 15 & 16 & 15 & 14 & 13 & 11 & 12 & 11 & 11 & 10 & 10 & 9\\
			\bottomrule
		\end{tabular}
		\caption{known results on the birational geometry of $\Mgn$. It is unirational (resp., uniruled) for $n\leq b(g)$ (resp., $n\leq \sigma(g)$); $\cM_{g,n}$ has non-negative (resp., maximal) Kodaira dimension for $n\geq \eta(g)$ (resp., $n\geq \tau(g)$).\label{tableBiratGeom}}
	\end{center}
\end{table}

We remark that in \cite{Logan} it was previously claimed that $\cM_{11,n}$ is unirational for $n\leq 10$. However, Barros in \cite{BarrosGeometry} noticed that the original argument contained a flaw and only proves the uniruledness of $\cM_{11,n}$ in that range. Barros managed to show that $\cM_{11,n}$ is unirational for $n \leq 6$, and that it is not unirational for $n=9$ or $n=10$. Thus, one cannot prove the unirationality of $\cM_{g,n}$ for $n$ up to $\sigma(g)$, as there certainly are cases which are uniruled but not unirational. However, it is not unreasonable to expect that the gap $\sigma(g)-b(g)$ in Table \ref{tableBiratGeom} can be reduced in many cases, especially for $g \geq 10$, where fewer results are known.

\subsection{New unirationality results for $\cM_{12,n}$ and $\cM_{13,n}$}

In this section we will prove the unirationality of $\cM_{12,n}$ for $n \leq 4$ and $\cM_{13,n}$ for $n \leq 3$. The strategy for proving these results is similar: we will exhibit a rational family of projective curves of genus $12$ (respectively, $13$) which is dominant on the corresponding moduli space $\cM_{12}$ (respectively, $\cM_{13})$. Both these families are constructed via liaison, respectively in $\PP^1 \times \PP^2$ and $\PP^3$; the construction allows us to impose a certain number of marked points on the corresponding curves.

Parts of the proofs are based on the explicit computation of single examples over a finite field, allowing us to show that some assumptions on the involved geometric objects, which correspond to open conditions, are generically satisfied.

\subsection{$\cM_{12,n}$}

The key step for proving the unirationality of $\cM_{12,4}$ will be the exploitation of a particular case of a construction by Gei\ss, which can be found in \cite{GeissThesis} and which we briefly recall in the next section.

\subsubsection{Gei\ss' construction}

In \cite{GeissThesis}, Gei\ss\ provided the proof of the unirationality of the Hurwitz spaces $\cH_{g,d}$, which parametrise $d$-sheeted simply branched covers of the projective line by smooth curves of genus $g$, up to isomorphism, for
\begin{itemize}[leftmargin=10pt]
    \item $d=6$ and $5 \leq g \leq 31$ or $g=33,34,35,36,39,40,45$;
\item $d=7$ and $6 \leq g \leq 12$.
\end{itemize}
In particular, several cases for $d=6$ and all the cases for $d=7$ were proved by establishing a correspondence between curves in $\PP^1\times \PP^2$ and certain submodules of the dual of their Hartshorne--Rao modules, in a similar fashion to Chang--Ran's approach in \cite{ChangRanUnirationality}.

We are mostly concerned with the case $(g,d)=(12,7)$, as we will use the dominant family exhibited by Gei\ss\ to deduce the unirationality of $\cM_{12,4}$. In what follows we briefly recall Gei\ss' construction for this specific case, in order to better present our argument and to provide some details which were omitted in \cite{GeissThesis}. The interested reader can find another specific case, namely $(g,d)=(10,6)$, examined in details in \cite[Appendix A]{CHGS}.

Let $C \rightarrow \PP^1$ be a general element in $\cH_{12,7}$. We consider a line bundle on $C$ of degree $10$ such that the map given by it and the assigned $g^1_7$ embeds $C$ in $\PP^1\times \PP^2$ as a curve of bidegree $(7,10)$. For general curves and line bundles, this embedding will be of maximal rank, as can be seen through the realisation of an explicit example. A Hilbert function computation shows that the truncated ideal $I':=(I_C)_{{}\geq (4,3)}$ admits a minimal free bigraded resolution over the Cox ring $R$ of $\PP^1\times \PP^2$ of the form
\[
0 \rightarrow R(-5,-5)^6 \xrightarrow{\varphi} R(-5,-4)^{10} \oplus R(-4,-5)^7 \rightarrow R(-4,-4)^9 \oplus R(-4,-3)^3 \rightarrow  I' \rightarrow 0.
\]
Let us denote by $F_i$ the terms in the above resolution, e.g., $F_1=R(-4,-4)^9 \oplus R(-4,-3)^3$. Then $K:=\coker \varphi^\vee$ is a module of finite length, called the \emph{truncated deficiency module}.

The first terms of a minimal free resolution of $K$ look like
\[
\dots \rightarrow G \xrightarrow{\psi} F_2^\vee  \xrightarrow{\varphi^\vee} F_3^\vee \rightarrow  K \rightarrow 0;
\]
one can prove that by composing $\psi$ with a general map $F_1^\vee \rightarrow G$ we obtain a matrix whose kernel is isomorphic to $R$. The entries of the corresponding induced map $R \rightarrow F_1^\vee$ generate $I'$, allowing us to recover $I'$ and thus the original curve $C$ from its truncated deficiency module $K$.

This correspondence can be exploited by constructing a rational family of modules $K$ and showing that such family leads to a family of curves $C$ in $\PP^1\times \PP^2$ which is dominant (by considering the first projection) on $\cH_{12,7}$. As it turns out, the main difficulty lies in the construction of $K$, as a general matrix $F_2^\vee \rightarrow F_3^\vee$ will produce a module of finite length which in general has Hilbert function different from the one $K$ is expected to have. In our specific case, $K$ will have dimension zero in all bidegrees but
\begin{equation}
\begin{gathered}
\label{realHilbertF}
\dim K_{(-5+i,-5)} = 6-i \quad \mbox{ for } 0 \leq i \leq 5,\\
\dim K_{(-5,-4)} = 8, \qquad \dim K_{(-4,-4)} = 4, \qquad \dim K_{(-5,-3)} = 6;
\end{gathered}
\end{equation}
the cokernel of a general matrix $R(5,4)^{10} \oplus R(4,5)^7 \rightarrow R(5,5)^6$, however, will be zero-dimensional in bidegree $(-4,-4)$. It is thus necessary to construct a rational family of matrices (i.e., a family of matrices parametrised by free parameters) whose cokernels have Hilbert function as in \eqref{realHilbertF}. Following Gei\ss' construction and by using Macaulay's inverse systems, we implemented in \cite{pointedCurvesCode} the construction of a family of curves of genus $12$ and bidegree $(7,10)$ in $\PP^1\times \PP^2$, dominant on $\cM_{12}$ and dense in an irreducible component $H_{12,(7,10)}$ of the Hilbert scheme of curves of that genus and degree in $\PP^1\times \PP^2$.

\subsubsection{The unirationality of $\cM_{12,4}$}

Let $C$ be a smooth curve of genus $g$ and bidegree $(d_1,d_2)$ in $\PP^1\times \PP^2$. Let $h_1, h_2$ be natural numbers and consider $\mathcal{O}_{C}(h_1,h_2)=\mathcal{O}_{C} \otimes (\mathcal{O}_{\PP^1}(h_1)\boxtimes \mathcal{O}_{\PP^2}(h_2))$. If $d_1h_1+d_2h_2>2g-2$, we have $\hhh^1(\mathcal{O}_{C}(h_1,h_2))=0$; in such case, by Riemann--Roch,
\[
\hhh^0(\mathcal{O}_{C}(h_1,h_2)) = d_1h_1+d_2h_2 +1 - g.
\]
We can thus compute the minimum (expected) number of hypersurfaces containing $C$ as $\max(0,\hhh^0(\mathcal{O}_{\PP^1\times \PP^2}(h_1,h_2))-\hhh^0(\mathcal{O}_{C}(h_1,h_2)))$, i.e.,
\begin{equation}
\label{expectedHyp}
\max(0,(h_1+1)(h_2+2)(h_2+1)/2-d_1h_1+d_2h_2 -1 +g).
\end{equation}

\begin{theorem}
	\label{m124}
	The moduli space $\cM_{12,n}$ is unirational for $n \leq 4$.
\begin{proof}
As explained in the previous section, we have a unirational component $H_{12,(7,10)}$ of the Hilbert scheme of curves of genus $12$ and bidegree $(7,10)$ in $\PP^1\times \PP^2$, dominant on $\cM_{12}$. By \eqref{expectedHyp}, a curve $C\in H_{12,(7,10)}$ lies on at least $2$ independent hypersurfaces of bidegree $(2,4)$; a concrete example will show that they are exactly $2$ for a general $C$. We can therefore consider the complete intersection curve cut out by two such hypersurfaces and link $C$ to a curve $C'$. If we assume that $C$ and $C'$ meet transversally and that $C'$ is smooth, by \eqref{liaisonP1P2} $C'$ will belong to a component $H_{4,(9,6)}$ of the Hilbert scheme of curves of genus $4$ and bidegree $(9,6)$ in $\PP^1\times \PP^2$.

We summarise the situation in the following diagram:
\begin{equation}
\label{firstdiagram}
\xymatrix{
\tilde{H}_{12,(7,10)} \ar@{-->}[r]^-\beta \ar@{-->}[d]^-\gamma&
H_{12,(7,10)} \ar@{-->}[r]^-\alpha &
\cM_{12}\\
\tilde{H}_{4,(9,6)} \ar@{-->}[r]^-\delta & H_{4,(9,6)}
},
\end{equation}
where the general element of $
\tilde{H}_{12,(7,10)}$ is a pair $(C,L)$ such that $C \in H_{12,(7,10)}$ is a curve and $L$ is a 2-dimensional subspace of $\HHH^0(\mathcal{I}_C(2,4))$, and analogously for $\tilde{H}_{4,(9,6)}$. The forgetful map $\alpha$ is dominant, $\beta$ is birational, the map $\gamma$ sends $(C,\langle Y_1,Y_2\rangle)$ to the geometrically linked curve $C'$ via the two hypersurfaces determined by $Y_1, Y_2$, and $\delta$ is the obvious forgetful morphism. By construction, every object in \eqref{firstdiagram} is unirational.

By \eqref{expectedHyp}, a general curve in $H_{4,(9,6)}$ will be contained in at least $6$ independent hypersurfaces of bidegree $(2,4)$ (and exactly $6$ in general, as shown by a concrete example). This means that we can perform liaison back: via two general hypersurfaces of bidegree $(2,4)$ containing it, we can link $C'$ to a curve $C''$, which will be again of genus $12$ and bidegree $(7,10)$ in general. By performing liaison forth and back we obtain again a unirational dominant family on $H_{12,(7,10)}$.

The fact that $C'$ is contained in at least $6$ hypersurfaces of bidegree $(2,4)$ allows us to impose the choice of (up to) $4$ general points on the so-constructed $C''$, obtained via liaison by choosing the hypersurfaces passing through those points. In mathematical terms, we can enrich diagram \eqref{firstdiagram} with the following one
\[
\xymatrix{
\left(H_{4,(9,6)}\times(\PP^1 \times \PP^2)^4\right)^{\sim}
\ar@{-->}[r]^-\zeta \ar@{-->}[d]^-\eta &
H_{4,(9,6)}\times(\PP^1 \times \PP^2)^4 \ar[r]^-\varepsilon & H_{4,(9,6)}\\
\cM_{12,4}
},
\]
where the general element of the first term is a $6$-tuple $(C',p_1,p_2,p_3,p_4,L')$ such that $C' \in H_{4,(9,6)}$ is a curve, the $p_i$'s are points in $\PP^1 \times \PP^2$ and $L'=\langle Y_1',Y_2'\rangle$ is a $2$-dimensional subspace of $\HHH^0(\mathcal{I}_{C'\cup (\cup_{i}p_i)}(2,4))$. The map $\varepsilon$ is the first projection and $\zeta$ is birational. The map $\eta$ sends such a $6$-tuple to $[C'',p_1,p_2,p_3,p_4]$, where $C''$ is the geometrically linked curve to $C'$ via the two hypersurfaces determined by $Y_1', Y_2'$. By construction, the $p_i$'s belong to $C''$, and any $4$-tuple of points on $C''$ arises this way. The composition of $\eta$ with the forgetful morphism $\cM_{12,4} \rightarrow \cM_{12}$ clearly factors through and dominates $H_{12,(7,10)}$, so $\eta$ is dominant and the unirationality of $\cM_{12,4}$ follows from the unirationality of the source of $\eta$.

The only things left to show are that, for a general curve $C$ in the family constructed above, the residual curve $C'$ obtained via liaison is indeed smooth and intersects $C$ transversally; and that the numbers of hypersurfaces of bidegree $(2,4)$ containing them is as expected. These are open conditions on the family, and can be checked through the realisation of a specific example, as we do in \cite{pointedCurvesCode}.
\end{proof}
\end{theorem}

\subsubsection{Another proof for the unirationality of $\cM_{12,3}$}
We present here another argument which yields the unirationality of $\cM_{12,n}$ for $n\leq 3$. Even if superseded by Theorem \ref{m124}, we believe it is of independent interest as it relies only on geometric arguments and does not depend on Gei\ss' construction and the homological approach used in it. As mentioned in Remark \ref{SernesiAlternative}, one could get a different proof by using a family of curves constructed by Sernesi in \cite{SernesiUnirationality}, obtained again via a homological approach.

Let $C$ be a general curve of genus $12$. Since $\rho(12,3,12)=0$, $C$ admits a $g^3_{12}$, which gives an embedding of $C$ as a curve of degree $12$ in $\PP^3$. Since $\cO_C(5)$ is non-special (i.e., $\hhh^1(\cO_C(5))=0$) by Riemann--Roch, $C$ is contained in at least $\hhh^0(\mathcal{I}_C(5))={\binom{5+3}{5}}-(5\cdot12+1-12)=7$ quintic hypersurfaces. Consider the complete intersection given by two general such hypersurfaces and suppose that the residual curve $C'$ is smooth and that $C$ and $C'$ intersect transversally; these are open conditions on the choice of $(C,\cO_C(1))\in \cW^3_{12,12}$. By \eqref{liaisonPr}, $C'$ is a curve of genus $15$ and degree $13$.

By Riemann--Roch, the Serre dual bundle $\omega_{C'}\otimes \cO_{C'}(-1)$ has a $5$-dimensional space of global sections and degree $15$. Hence, it is expected to embed the curve $C'$ in $\PP^4$ as a curve of degree $15$. In this embedding, $C'$ is contained in at least $4$ cubic hypersurfaces. Let $C^{\prime\prime}\subset \PP^4$ be the curve linked to $C'$ via the complete intersection of three such cubic hypersurfaces. By \eqref{liaisonPr}, $C^{\prime\prime}$ is a curve of genus $9$ and degree $12$. Again by Riemann--Roch we have $\hhh^0(\omega_{C^{\prime\prime}}\otimes \cO_{C^{\prime\prime}}(-1))=1$: this means that the Serre dual divisor corresponds to an element of $\cW^0_{9,4}$, or to the class of an effective divisor of degree $4$ on $C^{\prime\prime}$.

\begin{theorem}
	\label{M123}
	The moduli space $\cM_{12,n}$ is unirational for $n \leq 3$.
\begin{proof}
We claim that the above construction can be reversed, i.e, that there exists a chain of correspondences

\begin{equation}
\label{doubleLiaisonM123}
\xymatrix{
\cW^3_{12,12}
\ar@{<-->}^-L[r] &
\cW^3_{15,13}
\ar@{<-->}^-S[r] &
\cW^4_{15,15}
\ar@{<-->}^-L[r]&
\cW^4_{9,12}
\ar@{<-->}^-S[r] &
\cW^0_{9,4}
}
\end{equation}
These maps should be actually thought of as maps between some components of the spaces here above; maps labelled by $S$ correspond to considering the Serre dual model, while maps labelled by $L$ correspond to taking suitable linkages. The reversibility comes from the fact that all the open assumptions we made about the generality and smoothness of the objects involved and the transversality of the curves in liaison hold true in some specific examples, as verified in \cite{pointedCurvesCode}.

By Riemann--Roch, the known unirationality of $\cM_{9,4}$ implies the unirationality of $\cW^0_{9,4}$. Following \eqref{doubleLiaisonM123} backwards, we find a unirational family which is dominant on a component of $\cW^3_{12,12}$. Such component dominates $\cM_{12}$: to see this, it is sufficient to show that, in general, for the curves $C$ in $\cW^3_{12,12}$ obtained from the correspondences \eqref{doubleLiaisonM123} the Petri map $\HHH^0(\omega_C(-1)) \otimes \HHH^0(\mathcal{O}_C(1)) \rightarrow \HHH^0(\omega_C)$ is injective. This condition is open and has been checked for specific examples in \cite{pointedCurvesCode} by showing that there are no linear relations among the generators of $\Gamma_*(\omega_C)$ in degree $-1$.

To impose (up to) three extra points on a curve of genus $12$ we use the first map on the left in \eqref{doubleLiaisonM123}. We consider a curve $C'$ of genus $15$ and degree $13$ in $\PP^3$, which belongs to a unirational component $H_{15,13}$ of the corresponding Hilbert scheme of curves. As can be checked on a specific example, $C'$ lies on exactly 5 independent quintic hypersurfaces. Thus we have the following diagram
\[
\xymatrix{
\left(H_{15,13}\times(\PP^3)^3\right)^{\sim}
\ar@{-->}[r]^-\beta \ar@{-->}[d]^-\gamma &
H_{15,13}\times(\PP^3)^3 \ar[r]^-\alpha & H_{15,13}\\
\cM_{12,3}
},
\]
where the general element of the first term is a $5$-tuple $(C',p_1,p_2,p_3,L)$ such that $C' \in H_{15,13}$ is a curve, the $p_i$'s are points in $\PP^3$ and $L=\langle Y_1,Y_2\rangle$ is a $2$-dimensional subspace of $\HHH^0(\mathcal{I}_{C'\cup p_1 \cup p_2 \cup p_3}(5))$. The map $\alpha$ is the first projection, $\beta$ is birational and $\gamma$ is the map which sends a $5$-tuple as above to $[C,p_1,p_2,p_3]$, where $C$ is the curve geometrically linked to $C'$ via $Y_1,Y_2$. By construction, the $p_i$'s belong to $C$, and any triple of points on $C$ arises this way; the composition of $\gamma$ with the forgetful morphism $\cM_{12,3} \rightarrow \cM_{12}$ is dominant by the above discussion, hence $\gamma$ is dominant as well. The unirationality of $\cM_{12,3}$ then follows from the unirationality of the source of $\gamma$.
\end{proof}
\end{theorem}

\begin{remark}
	\label{SernesiAlternative}
	In \cite{SernesiUnirationality}, Sernesi constructed a rational family of curves of genus $12$ and degree $12$ in $\PP^3$ which is dominant on $\cM_{12}$, proving thus the unirationality of $\cM_{12}$. Another way to prove Theorem \ref{M123} is to use this family and perform a liaison with respect to $2$ quintics in $\PP^3$ forth and back using the first map in \eqref{doubleLiaisonM123}, as done in the proof of Theorem \ref{M123}.
\end{remark}

\subsection{$\cM_{13,n}$}
Let $C$ be a general curve of genus $13$. Since $\rho(13,3,13)=1$, $C$ admits a $g^3_{13}$, which gives an embedding of $C$ as a curve of degree $13$ in $\PP^3$. Since $\cO_C(5)$ is non-special by Riemann--Roch, $C$ is contained in at least $3={\binom{8}{5}}-(5\cdot13+1-13)$ quintic hypersurfaces. The complete intersection cut out by two such hypersurfaces links $C$ to a curve $C'$, which by \eqref{liaisonPr} has genus $10$ and degree $12$ provided that it is smooth and that intersects $C$ transversally.

Let $D'$ be the divisor associated to the embedding of $C'\subset \PP^3$. Riemann--Roch yields $\hhh^0(K_{C'}-D')=1$, and therefore the Serre dual divisor determines an element of $\cW^0_{10,6}$, whose general element can be interpreted as a curve $C'$ together with $6$ unordered points on it.

\begin{theorem}
	\label{M133}
	The moduli space $\cM_{13,n}$ is unirational for $n \leq 3$.
\begin{proof}
    The above construction can be reversed: indeed, a curve of genus $10$ and degree $12$ in $\PP^3$ lies on at least $5$ independent quintic hypersurfaces (and exactly $5$ for general curves, as showed by a concrete example). For the choice of $n\leq 3$ general points in $\PP^3$, there are at least $2$ such hypersurfaces passing through the $n$ points so that we can reverse the liaison construction to obtain a new genus $13$ curve together with up to $3$ marked points.
    
    More precisely, if we assume that $\cW^0_{10,6}$ is unirational, then by considering the Serre dual model we have a unirational family of curves which dominates a component $H_{10,12}$ of the Hilbert scheme of curves of genus $10$ and degree $12$ in $\PP^3$. We have the following diagram
    \[
    \xymatrix{
\left(H_{10,12}\times(\PP^3)^3\right)^{\sim}
\ar@{-->}[r]^-\beta \ar@{-->}[d]^-\gamma &
H_{10,12}\times(\PP^3)^3 \ar[r]^-\alpha & H_{10,12}\\
\cM_{13,3}
},
    \]
    where the general element of the first term is a $5$-tuple of the form $(C',p_1,p_2,p_3,L)$ such that $C'\in H_{10,12}$ is a curve, $p_i \in \PP^3$ for any $i$ and $L=\langle Y_1,Y_2\rangle$ is a $2$-dimensional subspace of $\HHH^0(\mathcal{I}_{C'\cup (\cup_{i}p_i)}(5))$. The map $\alpha$ is the first projection and $\beta$ is birational. The map $\gamma$ sends such a $5$-tuple to $[C,p_1,p_2,p_3]$, where $C$ is a genus $13$, degree $13$ curve, geometrically linked to $C'$ via the two hypersurfaces determined by $Y_1, Y_2$. By construction, the $p_i$'s belong to $C$, and any triple of points on $C$ arises this way. 
    
    If we prove that $\cW^0_{10,6}$ is unirational and that the unirational family provided by the above construction dominates $\cM_{13}$, we have that $\gamma$ dominates $\cM_{13,3}$ and the conclusion follows from the unirationality of the source of $\gamma$.
    
    A general element in $\cW^0_{10,6}$ gives $6$ unordered points on $C'$; conversely, by Riemann--Roch $6$ general points on a general curve of genus $10$ provide an element of $\cW^0_{10,6}$. On the one hand, $\cM_{10,6}$ is not known to be unirational (see Table \ref{tableBiratGeom}); on the other hand, we only need $6$ \emph{unordered} points, as we want to consider them as a divisor on the curve. Hence, the unirationality of $\cM_{10,6}^u$ suffices, a result already provided in \cite{BarrosGeometry} and for which we present an alternative proof in Remark \ref{alternativeBarros}.
    
    To prove that the so-constructed unirational family of curves $C$ of genus $13$ and degree $13$ dominates $\cM_{13}$, it is sufficient to show that in general the Petri map $\HHH^0(\omega_C(-1)) \otimes \HHH^0(\mathcal{O}_C(1)) \rightarrow \HHH^0(\omega_C)$ is injective. This open condition has been checked for one particular example in \cite{pointedCurvesCode} by showing that there are no linear relations among the generators of $\Gamma_*(\omega_C)$ in degree $-1$.
\end{proof}
\end{theorem}

\section{On the unirationality of $\cM_{g,n}^u$}
\label{uniratmgnusection}

In a similar fashion for $\cM_{g,n}$ in Section \ref{uniratM1213}, one can wonder for which pairs $(g,n)$ the space $\cM_{g,n}^u$ is unirational, uniruled, of non-negative Kodaira dimension, or of general type.
As it turns out by Riemann--Roch, for $n \geq g$ one has a rational dominant map $\cM_{g,n}^u \dashrightarrow \cP_{g,n}$, whose fibre over a point $(C,L)$ is $\PP(\HHH^0(L))$. Here $\cP_{g,n}$ denotes the universal Picard variety, parametrising smooth curves of genus $g$ together with a line bundle of degree $n$. This link between $\cM_{g,n}^u$ and $\cP_{g,n}$ can be, and in facts has been, used to characterise the birational geometry of $\cM_{g,n}^u$. In the following proposition we summarise all the relevant known contributions, which can be found in \cite{Logan, VerraUnirationality, FarkasPopa, BiniFontanariViviani, FarkasVerraUniversalJacobians, FarkasVerraUniversalThetaDivisor, CasalainaKassViviani, BarrosGeometry}.

\begin{proposition}
Let $g<22$ and $n \in \mathbb{N}$. Then
\begin{itemize}[leftmargin=10pt]
    \item $\cM_{g,n}^u$ is unirational for $g \leq 9$;
    \item $\cM_{10,6}^u$ and $\cM_{10,7}^u$ are unirational;
    \item $\cM_{g,n}^u$ is not unirational for $n \geq g \geq 10$;
    \item $\cM_{g,n}^u$ is uniruled for $n \geq g+1$, and for $g=10,11$ and $n < g$;
    \item $\kappa(\cM_{10,10}^u)=0$, $\kappa(\cM_{11,11}^u)=19$; $\kappa(\cM_{g,g}^u)=3g-3$ for g $\geq$ 12;
    \item $\cM_{g,n}^u$ is of general type for $g \geq 12$ and $f(g) \leq n \leq g-1$, where $f(g)$ is defined in Table \ref{farkasfg}.
    \end{itemize}
\begin{table}[h!bt]
	\begin{center}
		\setlength{\tabcolsep}{5pt}		
		\begin{tabular}{c|cccccccccc} \toprule
			$g$ & 12 & 13& 14& 15& 16& 17& 18& 19& 20& 21\\
			\midrule
			$f(g)$ & 10 & 11 & 10 & 10 & 9 & 9 & 9 & 7 & 6 & 4 \\
			\bottomrule
		\end{tabular}
		\caption{\label{farkasfg}}
	\end{center}
\end{table}
\end{proposition}

The unirationality of $\cM_{g,n}^u$ for $g\geq 10$ and small $n$, however, still remains mysterious, as no positive results are available besides the ones easily deducible from the unirationality of $\cM_{g,n}$, with the unique exception of $\cM_{10,6}^u$ and $\cM_{10,7}^u$.

\subsection{New unirationality results for $\cM_{g,n}^u$}

The aim of this section is to prove the following:

\begin{theorem}
	\label{mgnunordered}
	The moduli space $\cM_{g,n}^u$ is unirational for $(g,n)=(11,7)$,  $(12,5)$ or $(12,6)$.
\begin{proof}
Let us prove the unirationality of $\cM_{12,5}^u$ first. By considering $5$ general points $P$ on a general curve $C$ of genus $12$, the Serre dual divisor corresponding to $P$ gives an embedding of $C$ as a curve of degree $17$ in $\PP^6$. The number of quadrics of $\PP^6$ containing $C$ is at least $5$, so $C$ can be linked to a curve $C'$ via a complete intersection of type $(2^5)$. If we assume that $C'$ is smooth and meets $C$ transversally, $C'$ will have degree $15$ and genus $9$ by \eqref{liaisonPr}; moreover, $C'$ will be contained in at least $6$ quadrics. By \cite[Theorem 1.2]{VerraUnirationality}, the unique irreducible component of the Hilbert scheme of curves of genus $9$ and degree $15$ in $\PP^6$ which dominates $\cM_9$ is unirational. Since $\cW^6_{12,17}$ is irreducible, the conclusion then follows if we show that, for a general $C'$ in such component, the liaison works as expected and the curve $C$ obtained by reversing the liaison construction is non-degenerate. This condition and the assumptions we made on $C'$ above are open in moduli and can be checked through the realisation of particular examples, as we do in \cite{pointedCurvesCode}.

For the two remaining cases we adopt a common approach, an instance of which can be found in Example \ref{m126u} below; the general approach differs only for the numerology involved and can be explained as follows. Assume that $\cM_{g,m}^u$ is unirational and consider $m < g-3$ general points on a general curve $C$ of genus $g$, so that the Serre dual divisor of the $m$ points leads to a general element of $\cW^{g-m-1}_{g,2g-2-m}$. This dual divisor embeds $C$ as a curve of degree $d$ in a suitable projective space $\PP^r$, where $r=g-m-1$. We then search for possible reversible liaison constructions and try to impose the choice of further $m'$ points, in order to obtain a projective curve with $m'$ marked points on it. The projective model yields by construction $m$ unordered points, which together with the additional $m'$ provides the unirationality for $\cM_{g,m+m'}^u$.

In Table \ref{datamgnu} at the end of the section we collect the values for $g,m,d,r$, as well as the liaison type L we use, which will be always given by $r-1$ hypersurfaces of the same degree $h$. After choosing $h$, the genus $g'$ and degree $d'$ of the curve obtained via liaison are fixed. As it turns out, with these data we can impose only $m'=1$ further point on genus $g$ curves, hence $n=m+1$ in these cases.

Theorem \ref{mgnunordered} is proved as soon as we show that the assumptions we make on the occurring linkages (namely, the smoothness of $C'$ and the transversality of the intersections of $C$ and $C'$) are satisfied in general. These correspond to open conditions, which are then proved to be verified for general elements of the spaces under investigation via the construction of specific examples, done in \cite{pointedCurvesCode}.
\end{proof}
\end{theorem}

\begin{example}
\label{m126u}
For the unirationality of $\cM_{12,6}^u$ the aforementioned approach goes as follows. We can start from the unirationality of $\cM_{12,5}^u$ granted by Theorem \ref{mgnunordered}, so that $m=5$. A general element of $\cM_{12,5}^u$ provides a general curve $C$ of genus $g=12$ and degree $d=17$ in $\PP^{6}$. Such  curve is contained in at least $5$ quadric hypersurfaces, so we can consider $C'$ linked to $C$ via $(2^5)$; if we assume that $C'$ is smooth and meets $C$ transversally, then by \eqref{liaisonPr} $C'$ will have genus $g'=9$ and degree $d'=15$. This process gives a unirational family of curves dominating a component  $H_{9,15}$ of the Hilbert scheme of curves of genus $9$ and degree $15$ in $\PP^6$.

The curve $C'$ is contained in at least $6$ quadric hypersurfaces (and exactly $6$ in general, as shown by a concrete example), so that we can reverse the construction and impose $m'=1$ further point on a general curve $C''$ (which will be a priori different from $C$). More precisely, we have the following diagram
\[
    \xymatrix{
\left(H_{9,15}\times\PP^6\right)^{\sim}
\ar@{-->}[r]^-\beta &
H_{9,15}\times\PP^6 \ar[r]^-\alpha & H_{9,15}
},
    \]
    where the general element of the first term is a triple of the form $(C',p,L')$ such that $C'\in H_{9,15}$ is a curve, $p \in \PP^6$ and $L'=\langle Y_i\rangle_{i=1..5}$ is a $5$-dimensional subspace of $\HHH^0(\mathcal{I}_{C'\cup p}(2))$. The map $\alpha$ is the first projection, while $\beta$ is birational.
    
    To $(C',p,L')$ we can associate $[C'',p]$, where $C''$ is the geometrically linked curve to $C'$ via the hypersurfaces determined by $Y_i$. By construction, $p$ belongs to $C''$ and any point on $C''$ arises this way. We can define a rational map $\gamma$ as
    \[
    \gamma((C',p,L')) = (C'',\mathcal{O}_{C''}(p-H)\otimes \omega_{C''}),
    \]
    where $H$ is the hyperplane divisor corresponding to the embedding of $C''$ in $\PP^6$. By construction, $\gamma$ is dominant and $\mathcal{O}_{C''}(p-H)\otimes \omega_{C''}$ will in general be an effective divisor made up of $6$ unordered general points, so that $\gamma$ can be regarded as a rational map
    \[
    \gamma:\xymatrix{
    \left(H_{9,15}\times\PP^6\right)^{\sim} \ar@{-->}[r] & \cM_{12,6}^u.
    }
    \]
    The conclusion then follows from the unirationality of the source of $\gamma$.
\end{example}

\begin{remark}
\label{alternativeBarros}
We remark that the same general argument introduced in the proof of Theorem \ref{mgnunordered} and detailed in Example \ref{m126u} can be used to find an alternative, constructive proof of the unirationality of $\cM_{10,6}^u$ and $\cM_{10,7}^u$, which was already proved in \cite{BarrosGeometry}. The corresponding data for the liaison process are collected in Table \ref{datamgnu} as well.
\end{remark}

\begin{table}[h!tb]
	\begin{center}
		\begin{tabular}{c|cccccc|c} 
			\toprule
			$g$ &$m$ & $d$  &$r$ & L &$g'$&$d'$  &$n$ \\
			\midrule
			11&6&14 &4 &$3^3$ & $9$ & 13& 7\\
			12&5&17 &6 &$2^5$ & $9$ & 15& 6\\
			10&5&13 &4 &$3^3$ & 12 & 14& 6\\
			10&6&12 &3 &$5^2$ & 13 & 13& 7\\
			\bottomrule
		\end{tabular}
		\caption{data for liaison constructions for the unirationality of $\cM_{g,n}^u$, Theorem \ref{mgnunordered} and Remark \ref{alternativeBarros}.\label{datamgnu}}
	\end{center}
\end{table}

\frenchspacing



\end{document}